\documentclass[letterpaper, 10 pt, conference]{IEEEtran}

\usepackage{fancyhdr}
\usepackage{url}
\usepackage{amsmath} 
\usepackage{amssymb}  
\usepackage{amsthm}
\usepackage{lipsum,adjustbox}
\usepackage[font=footnotesize]{caption}
\usepackage{multirow}
\usepackage{cite}

\usepackage[shortlabels]{enumitem}
\renewcommand{\theenumi}{(\roman{enumi})}

\usepackage{subfig}

\DeclareMathOperator{\N}{\mathcal{N}}

\DeclareMathOperator{\F}{\mathcal{F}}
\DeclareMathOperator{\C}{\mathcal{C}}

\newcommand{\NN}{\mathbb{N}}

\newcommand{\real}{\mathbb{R}}
\newcommand{\realpos}{\mathbb{R}_{>0}}
\newcommand{\realnn}{\mathbb{R}_{\geq 0}}

\DeclareMathOperator{\DD}{\mathcal{D}}

\newtheorem{theorem}{Theorem}[section]
\newtheorem{proposition}[theorem]{Proposition}
\newtheorem{lemma}[theorem]{Lemma}
\newtheorem{corollary}[theorem]{Corollary}
\theoremstyle{remark}
\newtheorem{remark}{Remark}
\theoremstyle{definition}

\newtheorem{example}{Example}

\newcommand{\longthmtitle}[1]{\mbox{} \textit{(#1):}}
\newcommand{\until}[1]{\{1,\dots,#1\}}
\newcommand{\setdef}[2]{\{#1 \; | \; #2\}}

\newcommand{\map}[3]{#1:#2 \rightarrow #3}

\newcommand{\oprocendsymbol}{\hbox{$\bullet$}}
\newcommand{\oprocend}{\relax\ifmmode\else\unskip\hfill\fi\oprocendsymbol}

\allowdisplaybreaks

\title{Nesterov Acceleration for Equality-Constrained
  Convex Optimization
  via 
  Continuously Differentiable Penalty Functions
  \thanks{This work was supported by the ARPA-e
    Network Optimized Distributed Energy Systems (NODES) program
    through award DE-AR0000695 and NSF Award ECCS-1917177.}
}

\author{Priyank Srivastava \quad Jorge
  Cort\'{e}s
  \thanks{The authors are with the Department of Mechanical and
    Aerospace Engineering, UC San Diego,
    {\tt\small \{psrivast,cortes\}@ucsd.edu}}
}

\graphicspath{ {Figures/} }

\begin{document}

\maketitle
\thispagestyle{fancy}

\begin{abstract}
  We propose a framework to use Nesterov's accelerated method for
  constrained convex optimization problems. Our approach consists of
  first reformulating the original problem as an unconstrained
  optimization problem using a continuously differentiable exact
  penalty function. This reformulation is based on replacing the
  Lagrange multipliers in the augmented Lagrangian of the original
  problem by Lagrange multiplier functions. The expressions of these
  Lagrange multiplier functions, which depend upon the gradients of
  the objective function and the
  constraints, 
  can make the unconstrained penalty function non-convex in general
  even if the original problem is convex. We establish sufficient
  conditions on the objective function and the constraints of the
  original problem under which the unconstrained penalty function is
  convex. This enables us to use Nesterov's accelerated gradient
  method for unconstrained convex optimization and achieve a
  guaranteed rate of convergence which is better than the
  state-of-the-art first-order algorithms for constrained convex
  optimization. Simulations illustrate our results.
\end{abstract}

\section{Introduction}
Convex optimization problems arise in areas like signal processing,
control systems, estimation, communication, data analysis, and machine
learning. They are also useful to bound the optimal values of certain
nonlinear programming problems, and  to approximate their
optimizers. Due to their ubiquitous nature and importance, much effort
has been devoted to efficiently solve them.
This paper is motivated by the goal of designing fast 
methods that combine the simplicity and ease of gradient methods with
acceleration techniques to efficiently solve constrained optimization
problems.

\emph{Literature review:} Gradient descent is a widespread method to
solve unconstrained convex optimization problems.  However, gradient
descent suffers from slow convergence. To achieve local quadratic
convergence, one can use Newton's method~\cite{SB-LV:04}. Newton's
method uses second-order information of the objective function and
requires the inversion of the Hessian of the function.  In contrast,
the accelerated gradient descent method proposed by
Nesterov~\cite{YEN:83} uses only first-order information combined with
momentum terms~\cite{WS-SB-EJC:16,BS-SSD-WS-MIJ:19} to achieve an
optimal convergence rate. For constrained convex optimization,
generalizations of gradient algorithms include the projected gradient
descent~\cite{YN:18} (for simple set constraints where the projection
of any point can be computed in closed form) and (continuous-time)
saddle-point or primal-dual dynamics (for general constraints), see
e.g.,~\cite{TK:56,KA-LH-HU:58,AC-BG-JC:17-sicon,AC-EM-SHL-JC:18-tac}.
When the saddle function is strongly convex-strongly concave, the
primal-dual dynamics converges exponentially fast, see
e.g.,~\cite{JC-VKNL:12}. Recent
work~\cite{JC-SKN:19-jnls,GQ-NL:19,DD-MRJ:19,SSD-WH:19} has explored
the partial relaxation of the strong convexity requirement while
retaining the exponential convergence rate. A method with improved
rate of convergence for constrained problems is accelerated mirror
descent~\cite{WK-AB-PLB:15} which, however necessitates the choice of
an appropriate mirror map depending on the geometry of the problem and
requires that each update solves a constrained optimization problem
(which might be challenging itself).  Some
works~\cite{PEG-DPR:12,PA-RO:17,NKD-SZK-MRJ:17} have sought to
generalize Newton's method for equality constrained problems,
designing second-order updates that require the inversion of the
Hessian matrix of the augmented Lagrangian.  Similar to gradient
descent, a generalization of Nesterov's method for constrained convex
optimization described in~\cite{YN:18} uses the~projection for simple
set constraints.  Here we follow an alternative route involving
continuously differentiable exact penalty
functions~\cite{RF:70,TG-EP:79} to convert the original problem into
the unconstrained optimization of a nonlinear function. The
works~\cite{GdP-LG:89,SL:92,GDP:94} generalize these penalty functions
and establish, under appropriate assumptions on the constraint set,
complete equivalence between the solutions of the constrained and
unconstrained problems.  We employ these penalty functions to
reformulate the constrained convex optimization problem and identify
sufficient conditions under which the unconstrained problem is also
convex.  Our previous work~\cite{PS-JC:18-cdc} explores the
distributed computation of the gradient of the penalty function when
the objective is separable and the constraints are locally
expressible.

\emph{Statement of contributions:} We consider equality-constrained
convex optimization problems. Our starting point is the exact
reformulation of this problem as the optimization of an unconstrained
continuously differentiable function. We show via a counterexample
that the unconstrained penalty function might not be convex for any
value of the penalty parameter even if the original problem is
convex. This motivates our study of sufficient conditions on the
objective and constraint functions of the original problem for the
unconstrained penalty function to be convex.  Our results are based on
analyzing the positive semi-definiteness of the Hessian of the penalty
function.  We provide explicit bounds below which, for any value of
the penalty parameter, the penalty function is either convex or
strongly convex on the domain, resp.  Since the optimizers of
the unconstrained convex penalty function are the same as the
optimizers of the original problem, we deduce that the proposed
Nesterov implementation solves the original constrained problem with
an accelerated convergence rate starting from an arbitrary
  initial condition.  Finally, we establish that Nesterov's algorithm
  applied to the penalty function renders the feasible set forward
  invariant.  This, coupled with the fact that the penalty terms
  vanish on the feasible set, ensures that the accelerated convergence
  rate is also achieved from any feasible initialization.

\section{Preliminaries}\label{sec:prelim}
We collect here\footnote{Throughout the paper, we employ the following
  notation.  Let $\real$, $\realpos$, $\realnn$ and $\NN$ be the set of real,
  positive real, non-negative real and natural numbers, resp. 
   We use $\mathcal{X}^o$ to denote the interior of a set
  $\mathcal{X}$. $e_i^n$ denotes the $n-$dimensional unit vector in
  direction $i$.
  Given a matrix $A$, $\N(A)$ denotes its nullspace, $A^\top$ its
  transpose, $\|A \|$ its 2-norm, $\lambda_{\min}(A)$ and
  $\lambda_{\max}(A)$ its minimum and maximum eigenvalue,
  resp. If $A$ is positive semi-definite, we let
  $\lambda_2(A)$ denote the smallest positive eigenvalue, regardless
  of the multiplicity of the eigenvalue~$0$.
  Finally, $V^\perp$ denotes the orthogonal complement of the vector
  space~$V$.  } basic notions of convex analysis~\cite{RTR:70,SB-LV:04}
and optimization~\cite{DPB:99}.

\emph{Convex Analysis:} Let $\C \subseteq \real^n$ be a convex set.  A
function $f : \real^n \rightarrow \real $ is \emph{convex} on $\C$ if
$ f(\lambda x + (1-\lambda)y) \leq \lambda f(x) + (1-\lambda) f(y)$,
for all $ x,y \in \C$ and $\lambda \in [0,1]$.  Convex functions have
the property of having the same local and global minimizers.  A
continuously differentiable $f : \real^n \rightarrow \real $ is convex
on $\C$ \emph{iff} $f(y) \geq f(x)+(y-x)^\top \nabla f(x)$, for all $
x,y \in \C$. A twice differentiable function is convex \emph{iff} its
Hessian is positive semi-definite. A twice differentiable function $f:
\real^n \rightarrow \real$ is \emph{strongly convex} on $\C$ with
parameter $c \in \realpos$ \emph{iff} $ \nabla^2 f(x) \geq c I$ for
all $x \in \C$.

\emph{Constrained Optimization:} Consider the following nonlinear
optimization problem
\begin{equation}\label{eq:nl}
  \begin{aligned}
    & \min_{x \in \mathcal{D}}
    & & f(x) \\
    &\; \; \text{s.t.}  & & 
    h(x)=0,
  \end{aligned}
\end{equation}
where $f:\real^n \rightarrow \real, 
\; h: \real^n \rightarrow \real^p$ are twice continuously
differentiable functions with $p \leq n$ and $\mathcal{D} \subset
\real^n$ is a compact set which is regular (i.e., $\DD=
\overline{\mathcal{D}^o}$).  The feasible set of~\eqref{eq:nl} is $\F=
\setdef{x \in \mathcal{D}}{ h(x)=0}$.
\emph{Linear independence constraint qualification} (LICQ) holds at $x
\in \real^n$ if
$ \{\nabla h_k(x)\}_{k \in \until{p}}$ are linearly independent.

The Lagrangian $L:\real^n  \times \real^p \rightarrow
\real$ associated to~\eqref{eq:nl}~is 
\begin{align*}
  L(x,\mu) = f(x)+
  \mu^\top h(x) ,
\end{align*}
where 
$\mu \in \real^p$ is the Lagrange
multiplier (also called dual variable) associated with the
constraints. A Karush-Kuhn-Tucker (KKT)
point for~\eqref{eq:nl} is $(\bar{x},\bar{\mu})$ such that
\begin{align*}
  \nabla_xL(\bar{x},\bar{\mu})
  &=0 , \qquad
  \quad h(\bar{x})=0 .
\end{align*}
Under LICQ, the KKT conditions are necessary for a point to be locally
optimal.

\emph{Continuously Differentiable Exact Penalty Functions:} With exact
penalty functions, the idea is to replace the constrained optimization
problem~\eqref{eq:nl} by an equivalent unconstrained problem.  Here, 
we discuss continuously differentiable exact penalty functions
following~\cite{TG-EP:79,GdP-LG:89}.  The key observation is that one
can interpret a KKT tuple as establishing a relationship between a
primal optimal solution $\bar{x}$ and the dual optimal 
$\bar{\mu}$.  In turn, the following result introduces
multiplier functions that extend this relationship to any  $x\in
\real^n$.

\begin{proposition}\longthmtitle{Multiplier functions and their
    derivatives~\cite{GdP-LG:89}}\label{prop:lambda}
  Assume that LICQ is satisfied at all $x \in \mathcal{D}$. Define
  $N:\real^n \rightarrow \real^{p \times p}$ as $N(x) = \nabla
  h(x)^\top \nabla h(x)$.  Then $N(x)$ is a positive definite matrix
  for all $x \in \DD$. Given the function $x \mapsto \mu(x)$ defined
  by  $\mu(x)= -N^{-1}(x)  \nabla h(x)^\top \nabla f(x)$,
  the following holds
  \begin{enumerate}[(a)]
  \item if $(\bar{x},\bar{\mu})$ is a KKT point
    for~\eqref{eq:nl}, then 
    $\mu(\bar{x})=\bar{\mu}$;
  \item $\mu : \real^n \rightarrow \real^p $ is a continuously differentiable function.
  \end{enumerate}
\end{proposition}

The multiplier function can be used to replace the multiplier vector
in the augmented Lagrangian
to define a continuously differentiable exact penalty
function. Consider the
continuously differentiable function
$\map{f^\epsilon}{\real^n}{\real}$,
\begin{align}\label{eq:penalty}
  f^{\epsilon}(x) &=f(x)+ \mu(x)^\top h(x)
  + \frac{1}{\epsilon}\|h(x)\|^2 .
\end{align}
The next result shows when $f^\epsilon$ is an exact penalty function.

\begin{proposition}\longthmtitle{Continuously differentiable exact
    penalty function~\cite{GdP-LG:89}}\label{prop:exactness}
  Assume LICQ is satisfied at all $x \in \mathcal{D}$ and consider
  the unconstrained problem
  \begin{align}\label{eq:unc}
    \min_{x \in \mathcal{D}^o} f^\epsilon(x) .
  \end{align}
  Then, there exists $\bar{\epsilon}$ such that the set of global
    minimizers of~\eqref{eq:nl} and~\eqref{eq:unc} are equal for all
    $\epsilon \in (0,\bar{\epsilon}]$.
\end{proposition}

\section{Problem Statement}\label{sec:problem}
Consider the following convex optimization problem
\begin{equation}\label{eq:convex}
  \begin{aligned}
    & \min_{x \in \mathcal{D}}
    & & f(x)\\
    &\; \; \text{s.t.}
    & & Ax-b=0,
  \end{aligned}
\end{equation}
where $f:\real^n \rightarrow \real$ is a twice continuously
differentiable convex function and $\mathcal{D}$ is a convex set.
Here $A \in \real^{p \times n}$ and $b \in \real^p$ with $p <
n$. Without loss of generality, we assume $A$ has full row rank
(implying that LICQ holds at all~$x \in \real^n$).

Our aim is to design a Nesterov-like fast method to
solve~\eqref{eq:convex}.  We do this by reformulating the problem as
an unconstrained optimization using continuously differentiable
penalty function methods, cf.  Section~\ref{sec:prelim}.  Then, we
employ the Nesterov's accelerated gradient method to design
\begin{subequations}\label{eq:algorithm}
  \begin{align}
    x_{k+1}&=y_k - \alpha \nabla f^\epsilon(y_k),
    \label{eq:algorithm-a}
    \\
    a_{k+1}&=\frac{1+\sqrt{4a_k^2+1}}{2}, \label{algorithm-b}
    \\
    y_{k+1}&=x_{k+1}+\frac{a_k-1}{a_{k+1}}(x_{k+1}-x_k), \label{eq:algorithm-c}
  \end{align}
\end{subequations}
where $\alpha \in \realpos$ is the stepsize. If $f^\epsilon$ is convex
with Lipschitz gradient $L$ and the algorithm is initialized at an
arbitrary initial condition $x_0$ with $y_0=x_0$ and $a_0=1$, then
according to~\cite[Theorem 1]{YEN:83},
\begin{subequations}\label{eq:Nesterov}
\begin{align}\label{eq:Nesterov-convex}
  f^\epsilon(x_k) - f^\epsilon(x^*) \leq \dfrac{C}{(k+1)^2},
\end{align}
where $x^* \in \real^n$ is a global minimizer of $f^\epsilon$ and
$C \in \realnn$ is a constant dependant upon the initial condition and
$L$.  If $f^\epsilon$ is strongly convex with parameter
$s \in \realpos$, and~\eqref{algorithm-b} and~\eqref{eq:algorithm-c}
are replaced by
\begin{align}\tag{\ref{eq:algorithm}d}\label{eq:algorithm-d}
  y_{k+1}&=x_{k+1}+\frac{\sqrt{L}-\sqrt{s}}{\sqrt{L}+\sqrt{s}}(x_{k+1}-x_k),
\end{align}
then one has from~\cite[Theorem 2.2.1]{YN:18}
\begin{align}\label{eq:Nesterov-strongly-convex}
  f^\epsilon(x_k) - f^\epsilon(x^*) \leq C_s \exp \left(-k
    \sqrt{\frac{s}{L}} \right) ,
\end{align}
\end{subequations}
where $C_s \in \realnn$ is a constant dependant upon the initial
condition, $s$, and $L$.  The key technical point for this approach to
be successful is to ensure that the penalty function $f^\epsilon$ is
(strongly) convex.  Section~\ref{sec:convex} below shows that this is
indeed the case for suitable values of the penalty parameter under
appropriate assumptions on the objective and constraint functions of
the original problem~\eqref{eq:convex}.

\begin{remark}\longthmtitle{Distributed Algorithm
    Implementation}\label{re:distributed}
  We note here that the algorithm~\eqref{eq:algorithm} is amenable to
  distributed implementation if the objective function is separable
  and the constraints are locally coupled. In fact, our previous
  work~\cite{PS-JC:18-cdc} has shown how, in this case, the
  computation of the gradient of the penalty function
  in~\eqref{eq:algorithm-a} can be implemented in a distributed
  way. Based on this observation, one could use the framework proposed
  here for fast optimization of convex problems in a distributed way.
  To obtain fast convergence, one could also use second-order augmented
  Lagrangian methods, e.g.,~\cite{PA-RO:17,NKD-SZK-MRJ:17}, but their
  distributed implementation faces the challenge of computing the
  inverse of the Hessian of the augmented Lagrangian to update the
  primal and dual variables. Even if the Hessian is sparse for
  separable objective functions and local constraints, its inverse in
  general is not.
  \oprocend
\end{remark}

\section{Convexity of the Penalty Function}\label{sec:convex}
We start by showing that the continuously differentiable exact penalty
function $f^\epsilon$ defined in~\eqref{eq:penalty} might not be
convex even if the original problem~\eqref{eq:convex} is convex.  For
the convex problem~\eqref{eq:convex}, the penalty function takes the form
\begin{align}\label{eq:penalty_eq}
  f^\epsilon(x) & = f(x) 
  \\
  & \quad -([AA^\top]^{-1}A \nabla f(x) )^\top
  (Ax-b)+\frac{1}{\epsilon}\|Ax-b\|^2. \notag
\end{align}
A look at this expression makes it seem like a sufficiently small
choice of $\epsilon$ might make $f^\epsilon$ convex for all $x \in
\mathcal{D}$.  The following shows that this is always not the case.

\begin{example}\longthmtitle{Non-convex penalty
    function}\label{ex:example}
  Consider
  \begin{equation*}
    \begin{aligned}
      &\min\limits_{x \in \mathcal{D}}  & & x_1^4 + x_2^4 \\
      &\; \; \text{s.t.}  & & x_1+x_2=0.
    \end{aligned}
  \end{equation*}
  The optimizer is $(0,0)$.  The penalty function takes the form
  \begin{align*}
    f^\epsilon(x)=x_1^4 + x_2^4 + \mu(x)^\top (x_1+x_2) +
    \dfrac{1}{\epsilon} (x_1+x_2)^2,
  \end{align*}
  where $\mu(x)=-(2x_1^3+2x_2^3)$.
  The Hessian of this function is 
  \begin{align*}
    \nabla^2 f^\epsilon(x) \! = \! \!
    \left[
    \begin{matrix}
      -12x_1^2-12x_1x_2+\dfrac{2}{\epsilon} &
      -6x_1^2-6x_2^2+\dfrac{2}{\epsilon}
      \\
      -6x_1^2-6x_2^2+\dfrac{2}{\epsilon} &
      -12x_2^2-12x_1x_2+\dfrac{2}{\epsilon}
      \end{matrix}
      \right] \!.
    \end{align*}
    If $x_1=0$, then the determinant of $ \nabla^2 f^\epsilon(x)$
    evaluates to $-36x_2^4$, which is independent of
    $\epsilon$. Hence, $f^\epsilon$ cannot be made convex over any set
    containing the vertical axis. \oprocend
\end{example}

Example~\ref{ex:example} shows that the penalty function cannot always
be convexified by adjusting the value of~$\epsilon$. Intuitively, the
reason for this fact is that the term susceptible to be scaled in the
expression~\eqref{eq:penalty_eq} which depends on the
parameter~$\epsilon$ is not strongly convex. This implies that there
are certain subspaces where non-convexity arising from the term that
involve the Lagrange multiplier function cannot be countered. In turn,
these subspaces are defined by the kernel of the Hessian of the last
term in the expression~\eqref{eq:penalty_eq} of the penalty function.

These observations motivate our study of conditions on the objective
function and the constraints that guarantee that the penalty function
is convex. In our discussion, we start by providing sufficient
conditions for the convexity of the penalty function over $\DD$.

\subsection{Sufficient Conditions for Convexity over the Domain}
Here we provide conditions for the convexity of the penalty
function~$f^\epsilon$ by establishing the positive semi definiteness
of its Hessian. Throughout the section, we assume $f$ is three times
differentiable.  Note that the gradient and the Hessian of
$f^\epsilon$ are given, resp., by
\begin{subequations}
  \begin{align}
    \nabla f^\epsilon(x)&=\nabla f(x)- \nabla^2 f(x) A^\top
    [AA^\top]^{-1}(Ax-b) \notag
    \\
    & \quad - A^\top [AA^\top]^{-1} A \nabla f(x)
    +\frac{2}{\epsilon}A^\top(Ax-b).  \label{eq:gradient_eq}
    \\
    \nabla^2 f^\epsilon(x)& =\nabla^2 f(x) - W(x) - \nabla^2 f(x)A^\top
    [AA^\top]^{-1}A \notag
    \\
    & \quad - A^\top [AA^\top]^{-1}A \nabla^2 f(x) +
    \frac{2}{\epsilon}A^\top A, \label{eq:hessian_eq}
  \end{align}
\end{subequations}
where we use the short-hand notation
\begin{align}\label{eq:W}
  W(x)= \sum\limits_{i=1}^n \nabla_{x_i} \nabla^2 f(x) A^\top [A
  A^\top]^{-1}(Ax-b) e_i^{n^\top}.
\end{align}
The following result provides sufficient conditions under which the
penalty function~\eqref{eq:penalty_eq} is convex on $\DD$.

\begin{theorem}\longthmtitle{Convexity of the penalty function}\label{thm:eq}
  For the optimization problem~\eqref{eq:convex}, assume $ \nabla^2
  f(x) - W(x) \succ 0$ for all $x \in \DD$ and let
  \begin{align*}
    \bar{\epsilon} = \min\limits_{x \in \DD} \dfrac{ 2
      \lambda_{\min}(AA^\top ) \lambda_{\min}(\nabla^2f(x) - W(x)) }{
      \lambda^2_{\max} ( \nabla^2 f(x) ) +R(x)
      \lambda_{\min}(\nabla^2f(x) - W(x)) } ,
  \end{align*}
  where $R(x)= 2 \lambda_{\max} (\nabla^2 f(x)) - \lambda_{\min} (
  \nabla^2 f(x) - W(x) )$.  Then $f^\epsilon$ is convex on $\DD$ for
  all $\epsilon \in (0, \bar{\epsilon}]$ and consequently the
  convergence guarantee~\eqref{eq:Nesterov-convex} holds.
\end{theorem}
\begin{IEEEproof}
  For an arbitrary $x \in \DD$, we are interested in determining the
  conditions under which $\nabla^2 f^\epsilon(x) \succeq 0$, or in
  other words, $v^\top \nabla^2 f^\epsilon(x) v \geq 0$ for all
  $v \in \real^n$. From~\eqref{eq:hessian_eq},
  \begin{align}\label{eq:expression}
    v^\top \nabla^2 f^\epsilon(x) v
    &= \frac{2}{\epsilon} v^\top A^\top
      A v + v^\top (\nabla^2 f(x) -W(x))v
    \\
    & \quad - 2v^\top( \nabla^2
      f(x)A^\top [AA^\top]^{-1}A )v. \notag
  \end{align}
  Let us decompose $v$ as $v = v^{\|}+v^{\perp}$, where $v^{\|}$ is
  the component of $v$ in the nullspace $\N(A)$ of $A$ and $v^{\perp}$
  is the component orthogonal to it. Then~\eqref{eq:expression}
  becomes
  \begin{align*}
    v^\top \nabla^2 f^\epsilon(x) v
    &= \frac{2}{\epsilon} v^{\perp
      \top} A^\top A v^{\perp} + v^\top (\nabla^2 f(x) -W(x) )v
    \\
    & \quad -
      2v^{\| \top} \nabla^2 f(x)A^\top [A A^\top]^{-1} A v^{\perp}
    \\
    & \quad
    -2v^{\perp \top} \nabla^2 f(x)A^\top [A A^\top]^{-1} A v^{\perp} .
  \end{align*}
  Since $A^\top(AA^\top)^{-1}Av^\perp = v^\perp$, cf.~\cite[Theorem
  1.1.1]{SLC-CDM:09}, the above expression reduces to
  \begin{align*}
    & v^\top \nabla^2 f^\epsilon(x) v = \frac{2}{\epsilon} v^{\perp
      \top} A^\top A v^{\perp} + v^\top (\nabla^2 f(x) -W(x) )v
    \\
    & \quad - 2v^{\| \top} \nabla^2 f(x) v^{\perp} -2v^{\perp \top}
    \nabla^2 f(x) v^{\perp}
    \\
    & \geq \Big( \frac{2}{\epsilon} \lambda_2 (A^\top A ) - 2
    \lambda_{\max} (\nabla^2 f(x)) \Big) \| v^\perp \|^2
    \\
    &\quad + \lambda_{\min} ( \nabla^2 f(x) - W(x) ) ( \|v^\perp \|^2
    + \| v^\| \|^2 )
    \\
    & \quad - 2 \lambda_{\max} (\nabla^2 f(x)) \|v^\perp \| \|v^\| \|
    \\
    & = 
    \begin{bmatrix}
      \!  \|v^\perp \| \! \\ \! \| v^\| \| \!
    \end{bmatrix}^{\hspace*{-0.1cm} \top} \hspace*{-0.2cm} \underbrace{\begin{bmatrix} S(x) & \hspace*{-0.1cm}  - \lambda_{\max} (\nabla^2 f(x))  \\
        \!  - \lambda_{\max} (\nabla^2 f(x)) & \hspace*{-0.1cm}
        \lambda_{\min} ( \nabla^2 f(x) - W(x) ) \!
      \end{bmatrix}}_{P(x)} \hspace*{-0.1cm} \begin{bmatrix}
      \!  \|v^\perp \| \! \\ \! \| v^\| \| \!
    \end{bmatrix} \! \!,
  \end{align*}
  where $S(x)=\dfrac{2}{\epsilon} \lambda_{\min} (A A^\top ) -
  R(x)$. Therefore, we deduce that $\nabla^2 f^\epsilon (x) \succeq 0$
  if $\epsilon$ is such that $P(x) \succeq 0$. Being a
  $2\times2$-matrix, the latter holds if $S(x)$ and determinant of
  $P(x)$ are non-negative. The determinant is non-negative if and only
  if
  \begin{align*}
    \epsilon \leq \dfrac{2 \lambda_{\min} (AA^\top) \lambda_{\min} (
      \nabla^2 f(x) - W(x) )}{\lambda^2_{\max} (\nabla^2 f(x))+R(x)
      \lambda_{\min} ( \nabla^2 f(x) - W(x) )}.
  \end{align*}
  The above value of $\epsilon$ also ensures that $S(x) > 0$.
  Taking the minimum over all $x \in \DD$ completes the proof.
\end{IEEEproof}

\begin{remark}\longthmtitle{Differentiability of the objective
    function}
  Note that the implementation of~\eqref{eq:algorithm} requires the
  objective function $f$ to be twice continuously differentiable,
  while the definition of $W$ in~\eqref{eq:W} involves the third-order
  partial derivatives of $f$. We believe that an extension of
  Theorem~\ref{thm:eq} could be pursued in case the objective function
  is only twice differentiable using tools from nonsmooth analysis,
  e.g.,~\cite{FHC:83}, but we do not pursue it here for space
  reasons. \oprocend
\end{remark}
 
The next result provides sufficient conditions under which the penalty
function is strongly convex on $\DD$.

\begin{corollary}\longthmtitle{Strong convexity of the penalty
    function}\label{co:strong}
  For the optimization problem~\eqref{eq:convex}, assume $ \nabla^2 f(x) -
  W(x) \succeq cI$ for all $x \in \DD$ and let
  \begin{align*}
    \bar{\epsilon}_s =  \min\limits_{x \in \DD}\! \dfrac{
      2\lambda_{\min}(AA^\top ) (c-s) }{ \lambda^2_{\max}  ( \nabla^2
      f(x) ) \!+ \! 2 (c-s)\lambda_{\max} ( \nabla^2
      f(x)) \! - \!(c-s)^2\!  } .
  \end{align*}
  Then $f^\epsilon$ is strongly convex on $\DD$ with parameter $s \in
  (0,c)$ for all $\epsilon \in (0, \bar{\epsilon}_s]$ and the
  convergence guarantee~\eqref{eq:Nesterov-strongly-convex} holds.
\end{corollary}
\begin{IEEEproof}
  Let us decompose $\nabla^2 f(x) - W(x)$ as $\nabla^2 f(x) - W(x)=
  B(x)+ sI$. Since $\nabla^2 f(x) - W(x) \succeq cI$, it follows that
  $B(x) \succeq (c-s)I $. Establishing that the penalty function is
  strongly convex with parameter $s$ is equivalent to establishing
  that, for all $x \in \DD$, $v^\top (\nabla^2 f^\epsilon(x) - sI )
  v\geq 0$ for all $v \in \real^n$. Following the same steps as in the
  proof of Theorem~\ref{thm:eq}, one can verify that this is true if,
  for all $x \in \DD$, $\epsilon$ is less than or equal to
  \begin{align*}
    \dfrac{ 2\lambda_{\min}(AA^\top ) \lambda_{\min}(B(x)) }{
      \lambda^2_{\max} ( \nabla^2 f(x) ) \! + \! 2
      \lambda_{\min}(B(x))\lambda_{\max} ( \nabla^2 f(x)) \! - \!
      \lambda^2_{\min}(B(x)) } .
  \end{align*}
  Replacing $\lambda_{\min}(B(x))$ by $c-s$, it follows that the
  penalty function is strongly convex with parameter $s$ if $\epsilon
  \leq \bar{\epsilon}_s$.
\end{IEEEproof}

It is easy to verify that Example~\ref{ex:example} does not satisfy
the sufficient condition identified in Theorem~\ref{thm:eq}.  This
condition can be interpreted as requiring the original objective
function to be sufficiently convex to handle the non-convexity arising
from the penalty for being infeasible. Finding the value of
$\bar{\epsilon}$ still remains a difficult problem as computing
$\lambda_{\min}(\nabla^2 f(x) - W(x))$ for all $x \in \DD$ is not
straightforward.  The next result simplifies the conditions of
Theorem~\ref{thm:eq} for linear and quadratic programming problems.

\begin{corollary}\longthmtitle{Sufficient conditions for problems with
    linear and quadratic objective functions}\label{co:quadratic}
  \begin{enumerate}[(i)]
  \item If the objective function in problem~\eqref{eq:convex} is
    linear, then the penalty function is convex on $\real^n$ for all
    values of~$\epsilon$;
  \item If the objective function in problem~\eqref{eq:convex} is
    quadratic with Hessian $Q \succ 0$, then the penalty function is
    convex on $ \real^n$ for all $\epsilon \in (0, \bar{\epsilon}]$,
    where
    \begin{align*}
      \bar{\epsilon} = \dfrac{ 2\lambda_{\min}(AA^\top
        ) \lambda_{\min}(Q) }{ \lambda^2_{\max} ( Q ) + 2
        \lambda_{\min}(Q)\lambda_{\max} ( Q) - \lambda^2_{\min}(Q) }
      .        
    \end{align*}
  \end{enumerate}
  In either case, the convergence guarantee~\eqref{eq:Nesterov-convex}
  holds.
  \end{corollary}
\begin{IEEEproof}
  We present our arguments for each case separately. For case (i), we
  have $\nabla^2 f(x)=0$. Hence, 
  \begin{align*}
  \nabla^2 f^\epsilon(x)=\dfrac{2}{\epsilon}A^\top A,
  \end{align*}
    which means that $\nabla^2 f^\epsilon(x) \geq 0$ for all $x \in
    \real^n$.  For case (ii),
    \begin{align*}
      f(x)=\dfrac{1}{2}x^\top Q x + h^\top x,
    \end{align*}
    where $Q \in \real^{n \times n}$ and $h \in \real^n$. The expression for the Hessian of $f^\epsilon$ becomes
    \begin{align*}
      \nabla^2 f^\epsilon(x) &= Q + \frac{2}{\epsilon}A^\top A \notag
      - Q A^\top [AA^\top]^{-1} A
      \\
      & \quad - A^\top [AA^\top]^{-1} A Q.
    \end{align*}
    Clearly $W(x)=0$ for all $x \in \real^n$, and the result follows
    from Theorem~\ref{thm:eq}.
\end{IEEEproof}

Following Corollary~\ref{co:strong}, one can also state similar
conditions for the penalty function to be strongly convex in the case
of quadratic programs, but we omit them here for space reasons.  From
Corollary~\ref{co:quadratic}, ensuring that the penalty function
convex is easier when the objective function is quadratic. This
follows from the fact that $W(x)$, which depends on the third order
derivatives of the objection function, vanishes. Hence, in the
quadratic case, the condition in Theorem~\ref{thm:eq} requiring the
Hessian of the objective function to be greater than $W(x)$ for all
$x \in \DD$ is automatically satisfied. 
In what follows we provide a very simple approach for general
objective functions.

\subsection{Convexity over Feasible Set Coupled with Invariance}
Here we present a simplified version of the proposed approach, which
is based on the fact that inside the feasible set the values of the
penalty and the objective functions is the same. To build on this
observation, we start by characterizing the extent to which the
constraints are satisfied under the Nesterov's algorithm. 

\begin{lemma}\longthmtitle{Forward invariance of the feasible set
    under Nesterov's algorithm applied to the penalty function}\label{lemma:feasible}
  Consider the Nesterov's accelerated gradient
  algorithm~\eqref{eq:algorithm} applied to the penalty
  function~\eqref{eq:penalty_eq} for an arbitrary $\epsilon \ge
  0$. If the algorithm is initialized at $y_0 = x_0$, with $x_0$
  belonging to the feasible set $\F$, then $\{x_k\}_{k=0}^{\infty}$,
  $\{y_k\}_{k=0}^{\infty} \in \F$.
\end{lemma}
\begin{IEEEproof}
  We need to prove that $Ax_k=b$ and $Ay_k=b$ for all $k \geq 0$ if
  $Ax_0=Ay_0=b$. We use the technique of mathematical induction to
  prove this.  Since this clearly holds for $k=0$, we next prove that
  if $Ax_k=Ay_k=b$, then $Ax_{k+1}=Ay_{k+1}=b$.
  From~\eqref{eq:algorithm-a} and~\eqref{eq:gradient_eq}, we have
  \begin{align*}
    Ax_{k+1}=&Ay_k-\alpha A \nabla f^\epsilon (y_k)
    \\
    =&Ay_k \! -\! \alpha A (\nabla f(y_k) \! - \! \nabla^2 f(y_k) A^\top
    [AA^\top]^{-1}(Ay_k\! - \!b ) \notag
    \\
    &   -\! A^\top [AA^\top]^{-1} A
    \nabla f(y_k)
    \! + \! \frac{2}{\epsilon}A^\top(Ay_k \!- \!b)). \notag
  \end{align*}
  Substituting $Ay_k=b$, the above expression evaluates to $b$
  independent of $\epsilon \ge 0$. Then from~\eqref{eq:algorithm-c},
  one has $Ay_{k+1}=b$. Since the argument above is independent of the
  values of $a_k$ for all $k \in \NN$, it holds for the strongly
  convex case~\eqref{eq:algorithm-d} as well, thus completing the
  proof by induction.
\end{IEEEproof}

As a consequence of this result, if the trajectory starts in the
feasible set $\F$, then it remains in it forever. This observation
allows us to ensure the convergence rate guarantee for any convex
objective function.

\begin{corollary}\longthmtitle{Accelerated convergence with feasible
    initialization}
    For the optimization problem~\eqref{eq:convex} and arbitrary
    $\epsilon \ge 0$, the algorithm~\eqref{eq:algorithm} initialized
    in $\F$ enjoys the guarantee~\eqref{eq:Nesterov} on convergence to
    the optimal value.
\end{corollary}
\begin{IEEEproof}
  Note that $f^\epsilon(x) = f(x)$ whenever $Ax = b$, and hence by
  definition, is automatically (strongly) convex on $\F$ regardless of
  the value of $\epsilon$. The convergence guarantee follows from this
  fact together with Lemma~\ref{lemma:feasible}.
\end{IEEEproof}

\begin{remark}\longthmtitle{Robustness of the proposed approach}
  Given any $x_0 \in \real^n$, one can find a feasible initial point
  $x_0 -A^\top [AA^\top]^{-1}(Ax_0-b)$ by projecting $x_0$ onto the
  feasible set $\F$, and then implement Nesterov's accelerated method
  with the projected gradient as $(I - A^\top [AA^\top]^{-1} A)\nabla
  f(x)$. In fact, this projected gradient method coincides with the
  approach proposed here when evaluated over~$ \F$.  The advantage of
  our approach resides in the incorporation of error-correcting terms
  incorporating the value of $Ax-b$, cf.~\eqref{eq:gradient_eq}, that
  penalize any deviation from the feasible set and hence provide
  additional robustness in the face of disturbances.  By contrast, the
  projected gradient approach requires either an error-free execution
  or else, if error is present, the trajectory may leave and remain
  outside the feasible set unless repeated projections of the updated
  state are taken.  The inherent robustness property of the approach
  proposed here is especially important in the context of distributed
  implementations, cf. Remark~\ref{re:distributed}, where agents need
  to collectively estimate (and hence only implement approximations
  of) $A^\top [AA^\top]^{-1} A \nabla f(x)$ and taking the projection
  in a centralized fashion is not possible.  The approach proposed
  here can also be extended to problems with convex inequality
  constraints, cf.~\cite{GdP-LG:89}, whereas computing the projection
  in closed form is not possible for general convex constraints.
  \oprocend
\end{remark}

\begin{figure*}[t]
  \centering
  \subfloat[][]{\includegraphics[width=.45\linewidth]{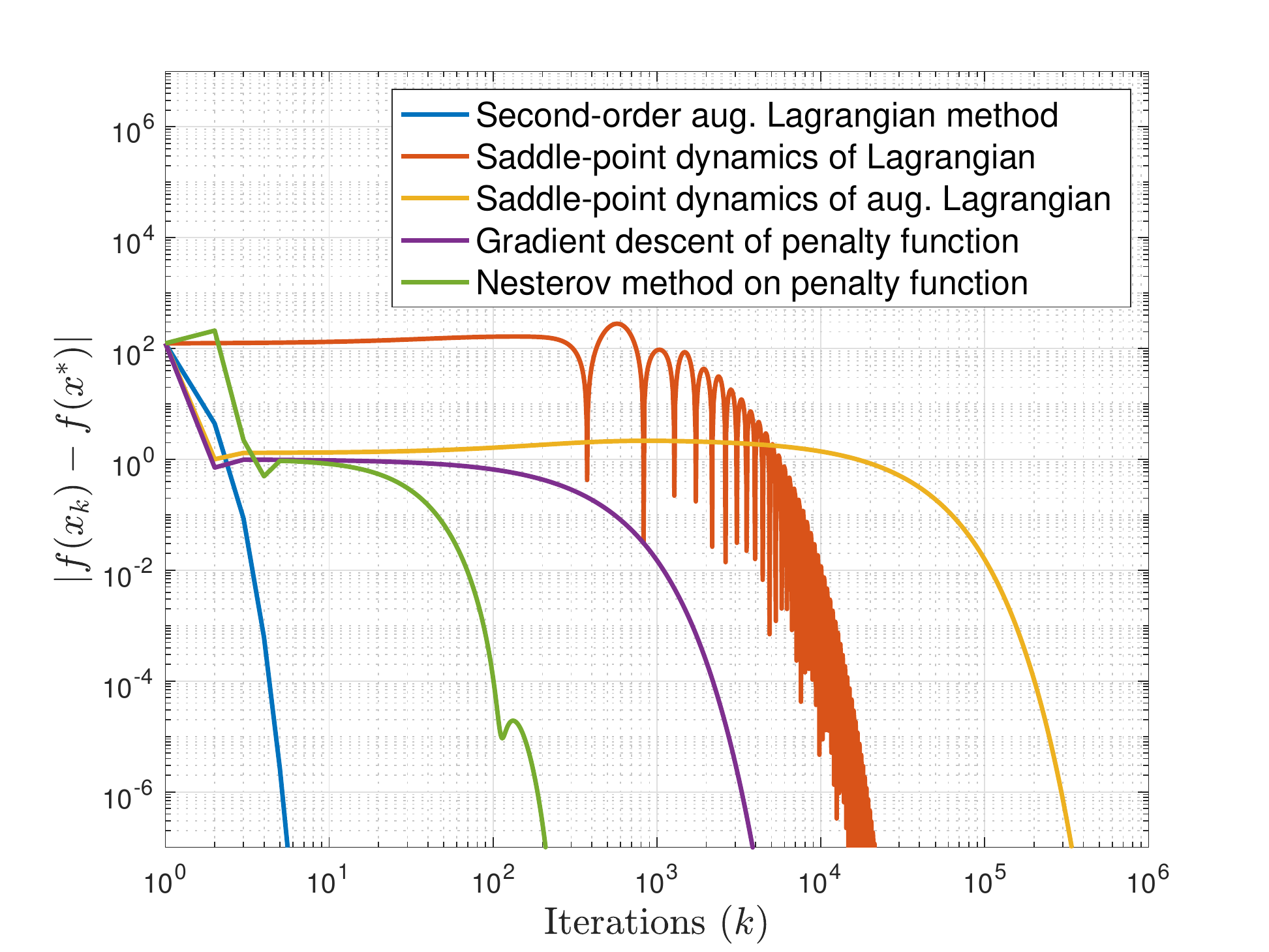}\label{fig:example_ra}}
  \quad
  \subfloat[][]{\includegraphics[width=.45\linewidth]{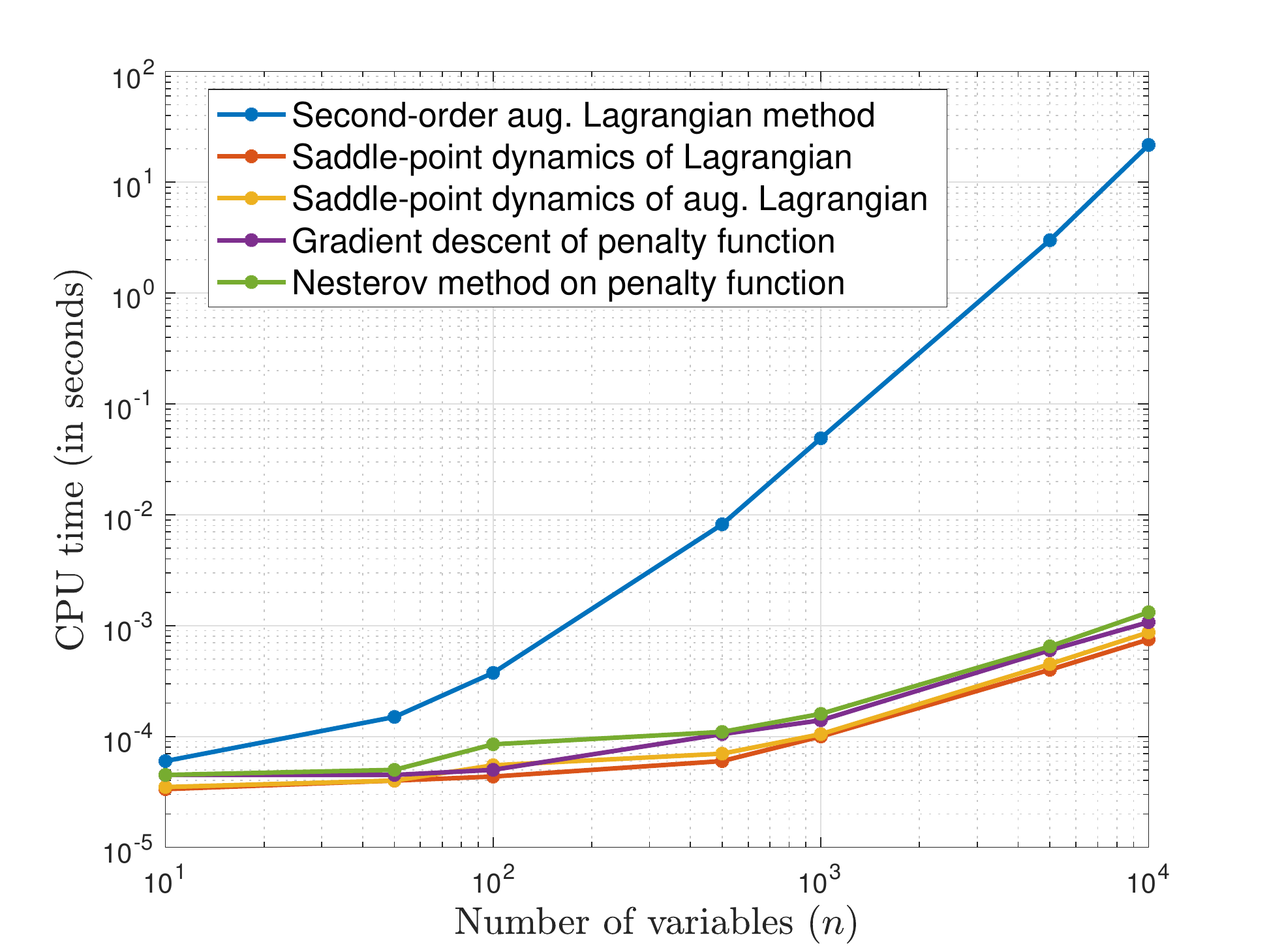}\label{fig:cpu_time}
  }
  \caption{Performance comparison of the proposed algorithm
    (Nesterov's acceleration on the penalty function) with the
    second-order augmented Lagrangian method~\cite{PA-RO:17}, the
    saddle-point dynamics~\cite{AC-EM-SHL-JC:18-tac,GQ-NL:19} applied
    to the Lagrangian and the augmented Lagrangian, respectively, and
    the gradient descent of the penalty function.  (a) shows the
    evolution of the error between the objective function and its
    optimal value for $n=50$ and (b) shows the computation time per
    iteration (note that the difference between second-order and
    first-order methods increases significantly with the problem
    dimension). For a desired level of accuracy, the proposed method
    outperforms the other methods when the number of iterations and
    the CPU time per iteration are jointly considered.}\label{fig:sim}
  \end{figure*}

\section{Simulations}\label{sec:sims}
In this section, we show the effectiveness of the proposed approach
through numerical simulations.  We consider 

\begin{equation*}
  \begin{aligned}
    & \min_{x \in \real^n} & & \sum\limits_{i=1}^n \frac{1}{2} \beta_i
    x^2_i + \gamma_i \exp(x_i)
    \\
    &\; \; \text{s.t.}  & & \sum\limits_{i=1}^n x_i=100,
  \end{aligned}
\end{equation*}
where $\beta_i$, $\gamma_i \in \realpos$.  We evaluate different
scenarios with values of $n$ as $10,50,100,500,1000, 5000$ and
$10000$.  We take $\DD=\setdef{x \in \real^n}{ \| x \|_{\infty} \leq
  5, \sum\limits_{i=1}^n x_i - 100 \leq 50 }$. By
Corollary~\ref{co:strong}, for $n=50$, the penalty function is
strongly convex on $\DD$ with parameter $s=0.01$ for all $\epsilon \in
(0, \bar{\epsilon}_s]$, where $\bar{\epsilon}_s=0.3603$. In our
simulations, we use $\epsilon = 10^{-1}$ and $\alpha = 10^{-3}$, resp.
Figure~\ref{fig:sim} compares the performance of the proposed method
with the second-order augmented Lagrangian method~\cite{PA-RO:17}, the
saddle-point dynamics~\cite{AC-EM-SHL-JC:18-tac,GQ-NL:19} applied to
the Lagrangian and the augmented Lagrangian, resp., and the gradient
descent applied to the penalty function.  Figure~\ref{fig:sim}(a)
shows the evolution of the error between the objective function and
its optimal value for $n=50$.  For the same level of accuracy, the
number of iterations taken by the second-order augmented Lagrangian
method is smaller by an order of magnitude compared to the proposed
method. However, one should note that the second-order augmented
Lagrangian method involves the inversion of Hessian, which becomes
increasingly expensive as the number of variables increases (see also
Remark~\ref{re:distributed}).  To illustrate this,
Figure~\ref{fig:sim}(b) shows the computation time per iteration of
the algorithms in Matlab version 2018a running on a Macbook Pro with
2GHz i5 processor and 8 GB ram. The time taken by the first-order
algorithms is about the same, and is smaller by several orders of
magnitude (depending on the number of variables) than the second-order
augmented Lagrangian method. When both aspects (number of iterations
and computation time per iteration) are considered together, the
proposed approach outperforms the other methods, especially if the
problem dimension is~large.

\section{Conclusions}\label{sec:conc}
We have presented a fast approach for constrained convex
optimization. We have provided sufficient conditions under which we
can reformulate the original problem as the unconstrained optimization
of a continuously differentiable convex penalty function. Our proposed
approach is based on the accelerated gradient method given by Nesterov
for unconstrained convex optimization, and has guaranteed convergence
rate when the penalty function is (strongly) convex. From simulations,
it is clear that in terms of computation time required to reach the
desired accuracy, the proposed method performs the best compared to
other state-of-the-art methods.  Based on our previous work, this
method is amenable to distributed optimization if, in the original
problem, the objective function is separable and the constraint
functions are locally expressible. Future work would explore the
effect of the choice of penalty parameter on the convergence speed of
the proposed strategy, the generalization of the conditions
identified here to ensure the penalty function is (strongly) convex
with inequality constraints, the extension of Nesterov's accelerated
gradient techniques to specific classes of non-convex functions (e.g.,
quasi-convex functions).


\end{document}